\documentclass[11pt,dvips,twoside,letterpaper]{article}
\usepackage[usenames]{color}
\usepackage{pslatex}
\usepackage{fancyhdr}
\usepackage{graphicx}
\usepackage{geometry}
\usepackage{amsmath}
\usepackage{amssymb}
\usepackage{amsfonts}
\usepackage{amsthm,amscd}

\def\figurename{Figure}
\makeatletter
\renewcommand{\fnum@figure}[1]{\figurename~\thefigure.}
\makeatother
\def\tablename{Table}
\makeatletter
\renewcommand{\fnum@table}[1]{\tablename~\thetable.}
\makeatother
\newtheorem{theorem}{Theorem}[section]
\newtheorem{lemma}[theorem]{Lemma}
\newtheorem{corollary}[theorem]{Corollary}

\theoremstyle{definition}
\newtheorem{definition}[theorem]{Definition}

\theoremstyle{remark}
\newtheorem{remark}[theorem]{Remark}
\numberwithin{equation}{section}

\begin{document}

\title{\bfseries\scshape{ A Hilbert-Type Integral Inequality in the Whole Plane Related to the
Hypergeometric Function and the Beta Function }}
\author{Michael Th. Rassias$^{1\ast }$ and Bicheng Yang$^{2}$ \and 1*.
Department of Mathematics, Princeton University,\and Fine Hall, Washington Road, Princeton, NJ 08544-1000, USA\and \& Department of Mathematics, ETH-Z\"{u}rich,\and R\"{a}mistrasse 101, 8092 Z\"{u}rich, Switzerland,\and michailrassias@math.princeton.edu \ michail.rassias@math.ethz.ch
\and 2. Department of Mathematics, Guangdong University of \and Education,
Guangzhou, Guangdong 510303, P. R. China \and E-mail: 1*. Corresponding
author: michael.rassias@math.ethz.ch \and 2. bcyang@gdei.edu.cn
bcyang818@163.com}

\date{}
\maketitle
\thispagestyle{empty} \setcounter{page}{1}

\renewcommand{\headrulewidth}{0pt}
\begin{abstract}
A new Hilbert-type integral inequality in the whole plane with the
non-homogeneous kernel and parameters is given. The constant factor related
to the hypergeometric function and the beta function is proved to be the
best possible. As applications, equivalent forms, the reverses, some
particular examples, two kinds of Hardy-type inequalities, and operator
expressions are considered.\newline


\textbf{\ }
\end{abstract}

\section{Introduction}

\footnotetext{\textbf{Foundation item:} This work is supported by the
National Natural Science Foundation of China (No. 61370186), and 2013
Knowledge Construction Special Foundation Item of Guangdong Institution of
Higher Learning College and University (No. 2013KJCX0140).} 

If $f(x),\:g(y)\geq 0,$ satisfy $$0<\int_{0}^{\infty }f^{2}(x)dx<\infty $$
and $$0<\int_{0}^{\infty }g^{2}(y)dy<\infty, $$ then we have the following
Hilbert's integral inequality (cf. \cite{HLP}):
\begin{equation}
\int_{0}^{\infty }\int_{0}^{\infty }\frac{f(x)g(y)}{x+y}dxdy<\pi \left(
\int_{0}^{\infty }f^{2}(x)dx\int_{0}^{\infty }g^{2}(y)dy\right) ^{\frac{1}{2}%
},  \label{1}
\end{equation}%
where the constant factor $\pi $ is the best possible. The inequality (\ref{1})
is very important in Mathematical Analysis and its applications (cf. \cite{HLP}, \cite{MPF}).
In recent years, by the use of the method of weight functions, a number of
extensions of (\ref{1}) were given by Yang (cf. \cite{BY1}). Noticing that
inequality (\ref{1}) is a homogenous kernel of degree -1, in 2009, a survey
of the study of Hilbert-type inequalities with the homogeneous kernels of
degree equal to negative numbers and some parameters is given in \cite{BY2}.
Recently, some inequalities with the homogenous kernels of degree 0 and
non-homogenous kernels have been proved (cf. \cite{BY3}-\cite{LB}). Other kinds of Hilbert-type inequalities are shown in \cite{MY1}-\cite%
{MY6}. All of the above integral inequalities are constructed in the quarter plane
of the first quadrant.

In 2007, Yang \cite{BY6} presented a new Hilbert-type integral inequality in the
whole plane, as follows:
\begin{eqnarray}
&&\int_{-\infty }^{\infty }\int_{-\infty }^{\infty }\frac{f(x)g(y)}{%
(1+e^{x+y})^{\lambda }}dxdy  \notag \\
&<&B(\frac{\lambda }{2},\frac{\lambda }{2})(\int_{-\infty }^{\infty
}e^{-\lambda x}f^{2}(x)dx\int_{-\infty }^{\infty }e^{-\lambda y}g^{2}(y)dy)^{%
\frac{1}{2}},  \label{2}
\end{eqnarray}%
where the constant factor $B(\frac{\lambda }{2},\frac{\lambda }{2})(\lambda
>0)$ is the best possible.\\ If $0<\lambda <1$, $p>1,$ $\frac{1}{p}+\frac{1}{q}%
=1,$ Yang \cite{BY8} derived another new Hilbert-type integral
inequality in the whole plane. Namely, he proved that
\begin{eqnarray}
&&\int_{-\infty }^{\infty }\int_{-\infty }^{\infty }\frac{1}{|1+xy|^{\lambda
}}f(x)g(y)dxdy  \notag \\
&<&k_{\lambda }\left[ \int_{-\infty }^{\infty }|x|^{p(1-\frac{\lambda }{2}%
)-1}f^{p}(x)dx\right] ^{\frac{1}{p}}\left[ \int_{-\infty }^{\infty }|y|^{q(1-%
\frac{\lambda }{2})-1}g^{q}(y)dy\right] ^{\frac{1}{q}},  \label{3}
\end{eqnarray}%
where the constant factor%
\begin{equation*}
k_{\lambda }=B(\frac{\lambda }{2},\frac{\lambda }{2})+2B(1-\lambda ,\frac{%
\lambda }{2})
\end{equation*}%
is still the best possible. Furthermore, Yang et al. \cite{BHY1}-\cite{ZZX2}
proved as well some new Hilbert-type integral inequalities in the whole plane.

In this paper, using methods from Real Analysis and by estimating the
weight functions, a new Hilbert-type integral inequality in the whole plane
with the non-homogeneous kernel and multi-parameters is shown, which gives an
extension of (\ref{3}). The constant factor related to the hypergeometric
function and the beta function is proved to be the best possible. As
applications, equivalent forms, the reverses, some particular examples,
two kinds of Hardy-type inequalities, and operator expressions are also
considered.

\section{Some Lemmas}
Initially, we introduce the following formula of the hypergeometric function $F$
(cf. \cite{WG}): If $\text{Re}(\gamma )>\text{Re}(\theta )>0,\:|\arg
(1-z)|<\pi ,$ $(1-zt)^{-\alpha }|_{z=0}=1,$ then
\begin{equation*}
F(\alpha ,\theta ,\gamma ,z):=\frac{\Gamma (\gamma )}{\Gamma (\theta )\Gamma
(\gamma -\theta )}\int_{0}^{1}t^{\theta -1}(1-t)^{\gamma -\theta
-1}(1-zt)^{-\alpha }dt,
\end{equation*}%
where, $$\Gamma (\eta )=\int_{0}^{\infty }x^{\eta -1}e^{-x}dx^{{}}(\text{Re}%
(\eta )>0)$$ is the gamma function. In particular, for $z=-1,\gamma
=\theta +1\:(\theta >0)$, $\alpha \in \mathbf{R},$ we have
\begin{equation}
\int_{0}^{1}t^{\theta -1}(1+t)^{-\alpha }dt=\frac{1}{\theta }F(\alpha
,\theta ,1+\theta ,-1)\in \mathbf{R}_{+}.  \label{4}
\end{equation}
\begin{lemma}
 If $\beta >-1,\mu ,\sigma >-\beta ,\mu +\sigma =\lambda
<1-\beta ,\delta \in \{-1,1\},$ we define two weight functions $\omega
(\sigma ,y)$ and $\varpi (\sigma ,x)$ as follows:
\begin{eqnarray}
\omega (\sigma ,y) &:&=\int_{-\infty }^{\infty }\frac{(\min \{1,|x^{\delta
}y|\})^{\beta }}{|1+x^{\delta }y|^{\lambda +\beta }}\frac{|y|^{\sigma }}{%
|x|^{1-\delta \sigma }}dx^{{}}(y\in \mathbf{R\backslash }\{0\}),  \label{5}
\\
\varpi (\sigma ,x) &:&=\int_{-\infty }^{\infty }\frac{(\min \{1,|x^{\delta
}y|\})^{\beta }}{|1+x^{\delta }y|^{\lambda +\beta }}\frac{|x|^{\delta \sigma
}}{|y|^{1-\sigma }}dy^{{}}(x\in \mathbf{R\backslash }\{0\}).  \label{6}
\end{eqnarray}%
Then we have
\begin{eqnarray}
\omega (\sigma ,y) &=&\varpi (\sigma ,x)=K(\sigma ):=\frac{1}{\beta +\sigma }%
F(\lambda +\beta ,\beta +\sigma ,1+\beta +\sigma ,-1)  \notag \\
&&+\frac{1}{\beta +\mu }F(\lambda +\beta ,\beta +\mu ,1+\beta +\mu ,-1)
\notag \\
&&+B(1-\lambda -\beta ,\beta +\sigma )+B(1-\lambda -\beta ,\beta +\mu )\in
\mathbf{R}_{+}.  \label{7}
\end{eqnarray}
\end{lemma}
\begin{proof}
 (i) For $\delta =1,$ by (\ref{5}) it follows
that%
\begin{equation*}
\omega (\sigma ,y)=\int_{-\infty }^{\infty }\frac{(\min \{1,|xy|\})^{\beta }%
}{|1+xy|^{\lambda +\beta }}\frac{|y|^{\sigma }}{|x|^{1-\sigma }}dx.
\end{equation*}

(a) If $y<0,$ setting $u=xy,$ we obtain%
\begin{eqnarray*}
\omega (\sigma ,y) &=&\int_{\infty }^{-\infty }\frac{(\min \{1,|u|\})^{\beta
}}{|1+u|^{\lambda +\beta }}\frac{(-y)^{\sigma }}{|u/y|^{1-\sigma }}\frac{1}{y%
}du \\
&=&\int_{-\infty }^{\infty }\frac{(\min \{1,|u|\})^{\beta }}{|1+u|^{\lambda
+\beta }}\frac{(-y)^{\sigma }(-y)^{1-\sigma }}{|u|^{1-\sigma }}\frac{1}{(-y)}%
du \\
&=&\int_{-\infty }^{\infty }\frac{(\min \{1,|u|\})^{\beta }|u|^{\sigma -1}}{%
|1+u|^{\lambda +\beta }}du;
\end{eqnarray*}

(b) if $y>0,$ setting $u=xy,$ it yields%
\begin{equation*}
\omega (\sigma ,y)=\int_{-\infty }^{\infty }\frac{(\min \{1,|u|\})^{\beta }}{%
|1+u|^{\lambda +\beta }}\frac{y^{\sigma }du}{|u/y|^{1-\sigma }y}%
=\int_{-\infty }^{\infty }\frac{(\min \{1,|u|\})^{\beta }|u|^{\sigma -1}}{%
|1+u|^{\lambda +\beta }}du.
\end{equation*}

(ii) For $\delta =-1,$ setting $X=x^{-1},$ we obtain%
\begin{eqnarray*}
\omega (\sigma ,y) &=&\int_{-\infty }^{0}\frac{(\min \{1,|x^{-1}y|\})^{\beta
}}{|1+x^{-1}y|^{\lambda +\beta }}\frac{|y|^{\sigma }}{|x|^{1+\sigma }}%
dx+\int_{0}^{\infty }\frac{(\min \{1,|x^{-1}y|\})^{\beta }}{%
|1+x^{-1}y|^{\lambda +\beta }}\frac{|y|^{\sigma }}{|x|^{1+\sigma }}dx \\
&=&\int_{0}^{-\infty }\frac{(\min \{1,|Xy|\})^{\beta }}{|1+Xy|^{\lambda
+\beta }}\frac{-|y|^{\sigma }dX}{|X^{-1}|^{1+\sigma }X^{2}}+\int_{\infty
}^{0}\frac{(\min \{1,|Xy|\})^{\beta }}{|1+Xy|^{\lambda +\beta }}\frac{%
-|y|^{\sigma }dX}{|X^{-1}|^{1+\sigma }X^{2}} \\
&=&\int_{-\infty }^{0}\frac{(\min \{1,|Xy|\})^{\beta }}{|1+Xy|^{\lambda
+\beta }}\frac{|y|^{\sigma }dX}{|X|^{1-\sigma }}+\int_{0}^{\infty }\frac{%
(\min \{1,|Xy|\})^{\beta }}{|1+Xy|^{\lambda +\beta }}\frac{|y|^{\sigma }dX}{%
|X|^{1-\sigma }} \\
&=&\int_{-\infty }^{\infty }\frac{(\min \{1,|xy|\})^{\beta }}{%
|1+xy|^{\lambda +\beta }}\frac{|y|^{\sigma }}{|x|^{1-\sigma }}dx.
\end{eqnarray*}
 Hence, for $\delta \in \{-1,1\},$ we obtain the following expression:%
\begin{eqnarray*}
\omega (\sigma ,y) &=&\int_{-\infty }^{\infty }\frac{(\min \{1,|u|\})^{\beta
}|u|^{\sigma -1}}{|1+u|^{\lambda +\beta }}du=K_{1}(\sigma )+K_{2}(\sigma ),
\\
K_{1}(\sigma ) &:&=\int_{-1}^{1}\frac{(\min \{1,|u|\})^{\beta }|u|^{\sigma
-1}}{(1+u)^{\lambda +\beta }}du, \\
K_{2}(\sigma ) &:&=\int_{\mathbf{R\backslash }(-1,1)}\frac{(\min
\{1,|u|\})^{\beta }|u|^{\sigma -1}}{|1+u|^{\lambda +\beta }}du.
\end{eqnarray*}

In view of (\ref{4}), and the following formula of the beta function (cf.
\cite{WG}):%
\begin{equation*}
B(p,q):=\int_{0}^{1}\frac{v^{p-1}}{(1-v)^{1-q}}dv^{{}}(p,q>0),
\end{equation*}%
we obtain%
\begin{eqnarray}
K_{1}(\sigma ) &=&\int_{-1}^{0}\frac{(-u)^{\beta +\sigma -1}}{(1+u)^{\lambda
+\beta }}du+\int_{0}^{1}\frac{u^{\beta +\sigma -1}}{(1+u)^{\lambda +\beta }}%
du  \notag \\
&=&\int_{0}^{1}\frac{v^{\beta +\sigma -1}}{(1-v)^{\lambda +\beta }}%
dv+\int_{0}^{1}\frac{u^{\beta +\sigma -1}}{(1+u)^{\lambda +\beta }}du  \notag
\\
&=&B(1-\lambda -\beta ,\beta +\sigma )+\frac{1}{\beta +\sigma }F(\lambda
+\beta ,\beta +\sigma ,1+\beta +\sigma ,-1).  \label{8}
\end{eqnarray}%
Setting $v=\frac{1}{u},$ it follows that%
\begin{eqnarray}
K_{2}(\sigma ) &=&\int_{\mathbf{-\infty }}^{-1}\frac{(\min
\{1,(-u)\})^{\beta }(-u)^{\sigma -1}}{|1+u|^{\lambda +\beta }}%
du+\int_{1}^{\infty }\frac{(\min \{1,u\})^{\beta }u^{\sigma -1}}{%
|1+u|^{\lambda +\beta }}du  \notag \\
&=&-\int_{0}^{-1}\frac{(\min \{1,(-\frac{1}{v})\})^{\beta }(-\frac{1}{v}%
)^{\sigma -1}}{|1+\frac{1}{v}|^{\lambda +\beta }v^{2}}dv-\int_{1}^{0}\frac{%
(\min \{1,\frac{1}{v}\})^{\beta }(\frac{1}{v})^{\sigma -1}}{|1+\frac{1}{v}%
|^{\lambda +\beta }v^{2}}dv  \notag \\
&=&\int_{-1}^{0}\frac{(\min \{1,(-v)\})^{\beta }(-v)^{\lambda -\sigma -1}}{%
|1+v|^{\lambda +\beta }}dv+\int_{0}^{1}\frac{(\min \{1,v\})^{\beta
}v^{\lambda -\sigma -1}}{|1+v|^{\lambda +\beta }}dv  \notag \\
&=&\int_{\mathbf{-1}}^{1}\frac{(\min \{1,|v|\})^{\beta }|v|^{\mu -1}}{%
(1+v)^{\lambda +\beta }}dv=K_{1}(\mu )  \notag \\
&=&B(1-\lambda -\beta ,\beta +\mu )+\frac{1}{\mu }F(\lambda +\beta ,\beta
+\mu ,1+\beta +\mu ,-1).  \label{9}
\end{eqnarray}

Setting $u=x^{\delta }y$ in (\ref{6}), for $x<0\:(x>0),$ we also get%
\begin{equation*}
\varpi (\sigma ,x)=\int_{-\infty }^{\infty }\frac{(\min \{1,|u|\})^{\beta
}|u|^{\sigma -1}}{|1+u|^{\lambda +\beta }}du=K(\sigma ).
\end{equation*}%
Hence we have (\ref{7}).
\end{proof}
\begin{remark}
 By Taylor's formula, we obtain
\begin{eqnarray*}
&&\frac{1}{\beta +\sigma }F(\lambda +\beta ,\beta +\sigma ,1+\beta +\sigma
,-1) \\
&=&\int_{0}^{1}\frac{u^{\beta +\sigma -1}du}{(1+u)^{\lambda +\beta }}%
=\int_{0}^{1}\sum_{k=0}^{\infty }\left( _{_{{}}k}^{-\lambda -\beta }\right)
u^{k+\beta +\sigma -1}du \\
&=&\int_{0}^{1}\sum_{k=0}^{\infty }(-1)^{k}\left( _{_{{}}k}^{\lambda +\beta
+k-1}\right) u^{k+\beta +\sigma -1}du \\
&=&\int_{0}^{1}\sum_{k=0}^{\infty }[\left( _{_{{}}2k}^{\lambda +\beta
+2k-1}\right) -\left( _{_{{}}2k+1}^{\lambda +\beta +2k}\right) u]u^{2k+\beta
+\sigma -1}du.
\end{eqnarray*}%
Since  we have%
\begin{equation*}
\left( _{_{{}}2k}^{\lambda +\beta +2k-1}\right) -\left(
_{_{{}}2k+1}^{\lambda +\beta +2k}\right) u=[1-\frac{(\lambda +\beta +2k)u}{%
2k+1}]\left( _{_{{}}2k}^{\lambda +\beta +2k-1}\right) ,
\end{equation*}%
there exists a large number $k_{0}\in \mathbf{N}_{0}=\mathbf{N}\cup \{0\},$
such that $\lambda +\beta +2k_{0}>0$, and for any $s\in \mathbf{N},$%
\begin{eqnarray*}
&&\left( _{_{{}}2(k_{0}+s)}^{\lambda +\beta +2(k_{0}+s)-1}\right) -\left(
_{_{{}}2(k_{0}+s+1)}^{\lambda +\beta +2(k_{0}+s)+1}\right) u \\
&=&[1-\frac{(\lambda +\beta +2k_{0}+2s)u}{2(k_{0}+s)+1}]\left(
_{_{{}}2(k_{0}+s)}^{\lambda +\beta +2(k_{0}+s)-1}\right)  \\
&=&[1-\frac{(\lambda +\beta +2k_{0}+2s)u}{2(k_{0}+s)+1}] \\
&&\times \frac{\lambda +\beta +2k_{0}+2s-1}{2k_{0}+2s}\cdots \frac{\lambda
+\beta +2k_{0}}{2k_{0}+1}\left( _{_{{}}2k_{0}}^{\lambda +\beta
+2k_{0}-1}\right) ,
\end{eqnarray*}%
\begin{eqnarray*}
&&1-\frac{(\lambda +\beta +2k_{0}+2s)u}{2(k_{0}+s)+1} \\
&\geq &1-\frac{\lambda +\beta +2k_{0}+2s}{2(k_{0}+s)+1}=\frac{1-\lambda
-\beta }{2(k_{0}+s)+1}>0(u\in (0,1]).
\end{eqnarray*}%
It follows that for any $s\in \mathbf{N},$ we have
\begin{equation*}
sgn(\left( _{_{{}}\,\,\,2(k_{0}+s)}^{\lambda +\beta +2(k_{0}+s)-1}\right)
-\left( _{_{{}}\,\,\,2(k_{0}+s+1)}^{\lambda +\beta +2(k_{0}+s)+1}\right)
u)=sgn\left( _{_{{}}\,\,\,2k_{0}}^{\lambda +\beta +2k_{0}-1}\right) .
\end{equation*}%
By Lebesgue's term by term integration theorem (cf. \cite{JK2}), we have%
\begin{eqnarray*}
&&\frac{1}{\sigma }F(\lambda +\beta ,\beta +\sigma ,1+\beta +\sigma ,-1) \\
&=&\sum_{k=0}^{\infty }\int_{0}^{1}[\left( _{_{{}}\,\,\,2k}^{\lambda +\beta
+2k-1}\right) -\left( _{_{{}}\,\,\,2k+1}^{\lambda +\beta +2k}\right) u]u^{2k+\beta
+\sigma -1}du \\
&=&\sum_{k=0}^{\infty }\left( _{_{{}}\,\,\,k}^{-\lambda -\beta }\right)
\int_{0}^{1}u^{k+\beta +\sigma -1}du=\sum_{k=0}^{\infty }\frac{1}{k+\beta
+\sigma }\left( _{_{{}}\,\,\,k}^{-\lambda -\beta }\right) .
\end{eqnarray*}

Similarly, since%
\begin{equation*}
(-1)^{k}\left( _{_{{}}\,\,\,\,\,\,k}^{-\lambda -\beta }\right) =(-1)^{k}\frac{(-\lambda
-\beta )(-\lambda -\beta -1)\cdots (-\lambda -\beta -k+1)}{k!}
\end{equation*}
\begin{equation*}
=(-1)^{2k}\frac{(\lambda +\beta )(\lambda +\beta +1)\cdots (\lambda +\beta
+k-1)}{k!}=\left( _{_{{}}\,\,\,\,\,\,k}^{\lambda +\beta +k-1}\right) ,
\end{equation*}%
we obtain%
\begin{eqnarray*}
&&B(1-\lambda -\beta ,\beta +\sigma )=\int_{0}^{1}\frac{u^{\beta +\sigma
-1}du}{(1-u)^{\lambda +\beta }}=\int_{0}^{1}\sum_{k=0}^{\infty
}(-1)^{k}\left( _{_{{}}\,\,\,k}^{-\lambda -\beta }\right) u^{k+\beta +\sigma -1}du
\\
&=&\int_{0}^{1}\sum_{k=0}^{\infty }\left( _{_{{}}\,\,\,k}^{\lambda +\beta
+k-1}\right) u^{k+\beta +\sigma -1}du.
\end{eqnarray*}%
There exists a large number $k_{1}\in \mathbf{N},$ such that $\lambda
+\beta +k_{1}>0.$\\ Hence, for any $s\in \mathbf{N},$ we have
\begin{equation*}
\left( _{_{{}}k_{1}+s}^{\lambda +\beta +k_{1}+s-1}\right) =\frac{\lambda
+\beta +k_{1}+s-1}{k_{1}+s}\cdots \frac{\lambda +\beta +k_{1}}{k_{1}+1}%
\left( _{_{{}}k_{1}}^{\lambda +\beta +k_{1}-1}\right) ,
\end{equation*}%
and then it follows that%
\begin{equation*}
sgn\left( _{_{{}}\,\,\,k_{1}+s}^{\lambda +\beta +k_{1}+s-1}\right)
=sgn\left( _{_{{}}\,\,\,k_{1}}^{\lambda +\beta +k_{1}-1}\right) .
\end{equation*}%
Still by Lebesgue's term by term integration theorem, we obtain%
\begin{equation*}
B(1-\lambda -\beta ,\beta +\sigma )=\sum_{k=0}^{\infty }\left(
_{_{{}}\,\,\,k}^{\lambda +\beta +k-1}\right) \int_{0}^{1}u^{k+\beta +\sigma
-1}du=\sum_{k=0}^{\infty }\frac{(-1)^{k}}{k+\beta +\sigma }\left(
_{_{{}}\,\,\,k}^{-\lambda -\beta }\right) .
\end{equation*}

Hence, we deduce the following series expressions:
\begin{eqnarray}
K_{1}(\sigma ) &=&2\sum_{k=0}^{\infty }\frac{1}{2k+\beta +\sigma }\left(
_{_{{}}\,\,\,2k}^{-\lambda -\beta }\right) ,  \label{10} \\
K_{2}(\sigma ) &=&2\sum_{k=0}^{\infty }\frac{1}{2k+\beta +\mu }\left(
_{_{{}}\,\,\,2k}^{-\lambda -\beta }\right) ,  \label{11} \\
K(\sigma ) &=&2\sum_{k=0}^{\infty }\frac{4k+2\beta +\lambda }{(2k+\beta
+\sigma )(2k+\beta +\mu )}\left( _{_{{}}\,\,\,2k}^{-\lambda -\beta }\right) .
\label{12}
\end{eqnarray}
\end{remark}
\begin{lemma}
 Suppose that $p>1,\frac{1}{p}+\frac{1}{q}=1,\beta >-1,\mu
,\sigma >-\beta ,\mu +\sigma =\lambda <1-\beta ,$ $\delta \in \{-1,1\},$ $%
K(\sigma )$ as indicated by (\ref{7}) (or (\ref{12})).\\ 
If $f(x)$ is a
non-negative measurable function in $\mathbf{R},$ then we have
\begin{eqnarray}
J &:&=\int_{-\infty }^{\infty }|y|^{p\sigma -1}\left[ \int_{-\infty
}^{\infty }\frac{(\min \{1,|x^{\delta }y|\})^{\beta }}{|1+x^{\delta
}y|^{\lambda +\beta }}f(x)dx\right] ^{p}dy  \notag \\
&\leq &K^{p}(\sigma )\int_{-\infty }^{\infty }|x|^{p(1-\delta \sigma
)-1}f^{p}(x)dx.  \label{13}
\end{eqnarray}
\end{lemma}
\begin{proof}
 By H\"{o}lder's inequality (cf. \cite{JK1})
and (\ref{5}), we derive that
\begin{eqnarray}
&&\left[ \int_{-\infty }^{\infty }\frac{(\min \{1,|x^{\delta }y|\})^{\beta }%
}{|1+x^{\delta }y|^{\lambda +\beta }}f(x)dx\right] ^{p}  \notag \\
&=&\left\{ \int_{-\infty }^{\infty }\frac{(\min \{1,|x^{\delta }y|\})^{\beta
}}{|1+x^{\delta }y|^{\lambda +\beta }}[\frac{|x|^{(1-\delta \sigma )/q}}{%
|y|^{(1-\sigma )/p}}f(x)][\frac{|y|^{(1-\sigma )/p}}{|x|^{(1-\delta \sigma
)/q}}]dx\right\} ^{p}  \notag \\
&\leq &\int_{-\infty }^{\infty }\frac{(\min \{1,|x^{\delta }y|\})^{\beta }}{%
|1+x^{\delta }y|^{\lambda +\beta }}\frac{|x|^{(1-\delta \sigma )(p-1)}}{%
|y|^{1-\sigma }}f^{p}(x)dx  \notag \\
&&\times \left[ \int_{-\infty }^{\infty }\frac{(\min \{1,|x^{\delta
}y|\})^{\beta }}{|1+x^{\delta }y|^{\lambda +\beta }}\frac{|y|^{(1-\sigma
)(q-1)}}{|x|^{1-\delta \sigma }}dx\right] ^{p-1}  \label{14}
\end{eqnarray}%
\begin{equation*}
=\frac{(\omega (\sigma ,y))^{p-1}}{|y|^{p\sigma -1}}\int_{-\infty }^{\infty }%
\frac{(\min \{1,|x^{\delta }y|\})^{\beta }}{|1+x^{\delta }y|^{\lambda +\beta
}}\frac{|x|^{(1-\delta \sigma )(p-1)}}{|y|^{1-\sigma }}f^{p}(x)dx.
\end{equation*}%
Then, by (\ref{7}) and Fubini's theorem (cf. \cite{JK2}), it follows that
\begin{eqnarray}
J &\leq &K^{p-1}(\sigma )\int_{-\infty }^{\infty }\left[ \int_{-\infty
}^{\infty }\frac{(\min \{1,|x^{\delta }y|\})^{\beta }}{|1+x^{\delta
}y|^{\lambda +\beta }}\frac{|x|^{(1-\delta \sigma )(p-1)}}{|y|^{1-\sigma }}%
f^{p}(x)dx\right] dy  \notag \\
&=&K^{p-1}(\sigma )\int_{-\infty }^{\infty }\varpi (\sigma
,x)|x|^{p(1-\delta \sigma )-1}f^{p}(x)dx.  \label{15}
\end{eqnarray}%
Hence, by (\ref{7}), inequality (\ref{13}) follows.
\end{proof}
\section{ Main Results and Applications}
\begin{theorem}
Let $p>1,\frac{1}{p}+\frac{1}{q}=1,\beta
>-1,\mu ,\sigma >-\beta ,$ $\mu +\sigma =\lambda <1-\beta ,$ $\delta \in
\{-1,1\},$ $K(\sigma )$ as indicated by (\ref{7}) (or (\ref{12})).\\ 
If $%
f,g\geq 0,$ satisfy $$0<\int_{-\infty }^{\infty }|x|^{p(1-\delta \sigma
)-1}f^{p}(x)dx<\infty $$ and $$0<\int_{-\infty }^{\infty }|y|^{q(1-\sigma
)-1}g^{q}(y)dy<\infty ,$$ then we have the following equivalent inequalities:
\begin{eqnarray}
I &:&=\int_{-\infty }^{\infty }\int_{-\infty }^{\infty }\frac{(\min
\{1,|x^{\delta }y|\})^{\beta }}{|1+x^{\delta }y|^{\lambda +\beta }}%
f(x)g(y)dxdy  \notag \\
&<&K(\sigma )\left[ \int_{-\infty }^{\infty }|x|^{p(1-\delta \sigma
)-1}f^{p}(x)dx\right] ^{\frac{1}{p}}\left[ \int_{-\infty }^{\infty
}|y|^{q(1-\sigma )-1}g^{q}(y)dy\right] ^{\frac{1}{q}},  \label{16}
\end{eqnarray}%
\begin{eqnarray}
J &:&=\int_{-\infty }^{\infty }|y|^{p\sigma -1}\left[ \int_{-\infty
}^{\infty }\frac{(\min \{1,|x^{\delta }y|\})^{\beta }}{|1+x^{\delta
}y|^{\lambda +\beta }}f(x)dx\right] ^{p}dy  \notag \\
&<&K^{p}(\sigma )\int_{-\infty }^{\infty }|x|^{p(1-\delta \sigma
)-1}f^{p}(x)dx.  \label{17}
\end{eqnarray}%
where, the constant factors $K(\sigma )$ and $K^{p}(\sigma )$ are the best
possible.

In particular, for $\delta =1$, we obtain the following equivalent inequalities:
\begin{eqnarray}
&&\int_{-\infty }^{\infty }\int_{-\infty }^{\infty }\frac{(\min
\{1,|xy|\})^{\beta }}{|1+xy|^{\lambda +\beta }}f(x)g(y)dxdy  \notag \\
&<&K(\sigma )\left[ \int_{-\infty }^{\infty }|x|^{p(1-\sigma )-1}f^{p}(x)dx%
\right] ^{\frac{1}{p}}\left[ \int_{-\infty }^{\infty }|y|^{q(1-\sigma
)-1}g^{q}(y)dy\right] ^{\frac{1}{q}},  \label{18}
\end{eqnarray}%
\begin{eqnarray}
&&\int_{-\infty }^{\infty }|y|^{p\sigma -1}\left[ \int_{-\infty
}^{\infty }\frac{(\min \{1,|xy|\})^{\beta }}{|1+xy|^{\lambda +\beta }}f(x)dx%
\right] ^{p}dy  \notag \\
&<&K^{p}(\sigma )\int_{-\infty }^{\infty }|x|^{p(1-\sigma )-1}f^{p}(x)dx.
\label{19}
\end{eqnarray}
\end{theorem}
\begin{proof}
If (\ref{14}) takes the form of equality for some $%
y\in (-\infty ,0)\cup (0,\infty )$, then, there exist constants $A$ and $B$, 
such that they are not all zero, and
\begin{equation*}
A\frac{|x|^{(1-\delta \sigma )(p-1)}}{|y|^{1-\sigma }}f^{p}(x)=B\frac{%
|y|^{(1-\sigma )(q-1)}}{|x|^{1-\delta \sigma }}\,^{{}}a.e.^{{}}\,in^{{}}\,%
\mathbf{R}.
\end{equation*}%
Let us suppose that $A\neq 0$ (otherwise $B=A=0$). Then it follows that
\begin{equation*}
|x|^{p(1-\delta \sigma )-1}f^{p}(x)=|y|^{q(1-\sigma )}\frac{B}{A|x|}%
^{{}}\,a.e.\,^{{}}in^{{}}\mathbf{R},
\end{equation*}%
which contradicts the fact that $$0<\int_{-\infty }^{\infty }|x|^{p(1-\delta
\sigma )-1}f^{p}(x)dx<\infty.$$ Hence (\ref{14}) takes the form of strict
inequality. So does (\ref{15}), and we obtain (\ref{17}).

By H\"older's inequality (cf. \cite{JK1}), we have
\begin{eqnarray}
I &=&\int_{-\infty }^{\infty }\left[ |y|^{\sigma -\frac{1}{p}}\int_{-\infty
}^{\infty }\frac{(\min \{1,|x^{\delta }y|\})^{\beta }}{|1+x^{\delta
}y|^{\lambda +\beta }}f(x)dx\right] (|y|^{\frac{1}{p}-\sigma }g(y))dy  \notag
\\
&\leq &J^{\frac{1}{p}}\left[ \int_{-\infty }^{\infty }|y|^{q(1-\sigma
)-1}g^{q}(y)dy\right] ^{\frac{1}{q}}.  \label{20}
\end{eqnarray}%
Then by (\ref{17}), we get (\ref{16}). On the other hand, suppose that (\ref%
{16}) is valid. We then set
\begin{equation}
g(y):=|y|^{p\sigma -1}\left[ \int_{-\infty }^{\infty }\frac{(\min
\{1,|x^{\delta }y|\})^{\beta }}{|1+x^{\delta }y|^{\lambda +\beta }}f(x)dx%
\right] ^{p-1}\,\,\,(y\neq 0),  \label{21}
\end{equation}%
and then $$J=\int_{-\infty }^{\infty }|y|^{q(1-\sigma )-1}g^{q}(y)dy.$$ 
By (%
\ref{15}), we have $J<\infty $. If $J=0$, then (\ref{17}) is trivially
true; if $0<J<\infty $, then by (\ref{16}), we obtain
\begin{eqnarray}
0 &<&\int_{-\infty }^{\infty }|y|^{q(1-\sigma )-1}g^{q}(y)dy=J=I  \notag \\
&<&K(\sigma )\left[ \int_{-\infty }^{\infty }|x|^{p(1-\delta \sigma
)-1}f^{p}(x)dx\right] ^{\frac{1}{p}}\left[ \int_{-\infty }^{\infty
}|y|^{q(1-\sigma )-1}g^{q}(y)dy\right] ^{\frac{1}{q}}<\infty ,  \label{22}
\end{eqnarray}%
\begin{equation}
J^{\frac{1}{p}}=\left[ \int_{-\infty }^{\infty }|y|^{q(1-\sigma
)-1}g^{q}(y)dy\right] ^{\frac{1}{p}}<K(\sigma )\left[ \int_{-\infty
}^{\infty }|x|^{p(1-\delta \sigma )-1}f^{p}(x)dx\right] ^{\frac{1}{p}}.
\label{23}
\end{equation}%
Hence, we obtain (\ref{17}), which is equivalent to (\ref{16}).

We set $E_{\delta }:=\{x\in \mathbf{R};|x|^{\delta }\geq 1\},$ and $%
E_{\delta }^{+}:=E_{\delta }\cap \mathbf{R}_{+}=\{x\in \mathbf{R}%
_{+};x^{\delta }\geq 1\}.$ For $\varepsilon >0,$ we define two functions $\tilde{%
f}(x),\tilde{g}(y)$ as follows:
\begin{eqnarray*}
\tilde{f}(x) &:&=\left\{
\begin{array}{ll}
|x|^{\delta (\sigma -\frac{2\varepsilon }{p})-1}, & x\in E_{\delta } \\
0, & x\in \mathbf{R}\backslash E_{\delta }%
\end{array}%
\right. , \\
\tilde{g}(y) &:&=\left\{
\begin{array}{ll}
0, & y\in (-\infty ,-1]\cup \lbrack 1,\infty ) \\
|y|^{\sigma +\frac{2\varepsilon }{q}-1}, & y\in (-1,1)%
\end{array}%
\right. .
\end{eqnarray*}%
Then we obtain
\begin{eqnarray*}
\tilde{L} &:&=\left[ \int_{-\infty }^{\infty }|x|^{p(1-\delta \sigma )-1}%
\tilde{f}^{p}(x)dx\right] ^{\frac{1}{p}}\left[ \int_{-\infty }^{\infty
}|y|^{q(1-\sigma )-1}\tilde{g}^{q}(y)dy\right] ^{\frac{1}{q}} \\
&=&2\left( \int_{E_{\delta }^{+}}x^{-2\delta \varepsilon -1}dx\right) ^{%
\frac{1}{p}}\left( \int_{0}^{1}y^{2\varepsilon -1}dy\right) ^{\frac{1}{q}}=%
\frac{1}{\varepsilon }.
\end{eqnarray*}

We find%
\begin{eqnarray*}
I(x) &:&=\int_{-1}^{1}\frac{(\max \{1,|x^{\delta }y|\})^{\beta }}{%
|1+x^{\delta }y|^{\lambda +\beta }}|y|^{\sigma +\frac{2\varepsilon }{q}-1}dy
\\
&&\overset{y=-Y}{=}\int_{-1}^{1}\frac{(\max \{1,|-x^{\delta }Y|\})^{\beta }}{%
|1-x^{\delta }Y|^{\lambda +\beta }}|-Y|^{\sigma +\frac{2\varepsilon }{q}%
-1}dY=I(-x),
\end{eqnarray*}%
and then $I(x)$ is an even function. It follows that
\begin{eqnarray*}
&&\tilde{I}:=\int_{-\infty }^{\infty }\int_{-\infty }^{\infty }\frac{(\min
\{1,|x^{\delta }y|\})^{\beta }}{|1+x^{\delta }y|^{\lambda +\beta }}\tilde{f}%
(x)\tilde{g}(y)dxdy \\
&=&\int_{E_{\delta }}|x|^{\delta (\sigma -\frac{2\varepsilon }{p}%
)-1}I(x)dx=2\int_{E_{\delta }^{+}}x^{\delta (\sigma -\frac{2\varepsilon }{p}%
)-1}I(x)dx \\
&&\overset{u=x^{\delta }y}{=}2\int_{E_{\delta }^{+}}x^{-2\delta \varepsilon
-1}\left[ \int_{-x^{\delta }}^{x^{\delta }}\frac{(\min \{1,|u|\})^{\beta }}{%
|1+u|^{\lambda +\beta }}|u|^{\sigma +\frac{2\varepsilon }{q}-1}du\right] dx.
\end{eqnarray*}%
Setting $v=x^{\delta }$ in the above integral, by Fubini's theorem (cf. \cite%
{JK2}), we find%
\begin{eqnarray*}
&&\tilde{I}=2\int_{1}^{\infty }v^{-2\varepsilon -1}\left[ \int_{-v}^{v}\frac{%
(\min \{1,|u|\})^{\beta }}{|1+u|^{\lambda +\beta }}|u|^{\sigma +\frac{%
2\varepsilon }{q}-1}du\right] dv \\
&=&2\int_{1}^{\infty }v^{-2\varepsilon -1}\left\{ \int_{0}^{v}[\frac{(\min
\{1,u\})^{\beta }}{|1-u|^{\lambda +\beta }}+\frac{(\min \{1,u\})^{\beta }}{%
(1+u)^{\lambda +\beta }}]u^{\sigma +\frac{2\varepsilon }{q}-1}du\right\} dv
\\
&=&2\int_{1}^{\infty }v^{-2\varepsilon -1}\left\{ \int_{0}^{1}[\frac{1}{%
(1-u)^{\lambda +\beta }}+\frac{1}{(1+u)^{\lambda +\beta }}]u^{\beta +\sigma +%
\frac{2\varepsilon }{q}-1}du\right\} dv \\
&&+2\int_{1}^{\infty }v^{-2\varepsilon -1}\left\{ \int_{1}^{v}[\frac{1}{%
(u-1)^{\lambda +\beta }}+\frac{1}{(1+u)^{\lambda +\beta }}]u^{\sigma +\frac{%
2\varepsilon }{q}-1}du\right\} dv\\
&=&\frac{1}{\varepsilon }\int_{0}^{1}[\frac{1}{(1-u)^{\lambda +\beta }}+%
\frac{1}{(1+u)^{\lambda +\beta }}]u^{\beta +\sigma +\frac{2\varepsilon }{q}%
-1}du \\
&&+2\int_{1}^{\infty }(\int_{u}^{\infty }v^{-2\varepsilon -1}dv)[\frac{1}{%
(u-1)^{\lambda +\beta }}+\frac{1}{(1+u)^{\lambda +\beta }}]u^{\sigma +\frac{%
2\varepsilon }{q}-1}du \\
&=&\frac{1}{\varepsilon }\left\{ \int_{0}^{1}[\frac{1}{(1-u)^{\lambda +\beta
}}+\frac{1}{(1+u)^{\lambda +\beta }}]u^{\beta +\sigma +\frac{2\varepsilon }{q%
}-1}du\right. \\
&&\left. +\int_{1}^{\infty }[\frac{1}{(u-1)^{\lambda +\beta }}+\frac{1}{%
(1+u)^{\lambda +\beta }}]u^{\sigma -\frac{2\varepsilon }{p}-1}du\right\} .
\end{eqnarray*}

If the constant factor $K(\sigma )$ in (\ref{16}) is not the best possible,
then, there exists a positive number $k,$ with $K(\sigma )<k$, such that (%
\ref{16}) is valid when replacing $K(\sigma )$ by $k$. Then in particular,
we have $\varepsilon \tilde{I}<\varepsilon k\tilde{L},$ and
\begin{eqnarray}
&&\int_{0}^{1}[\frac{1}{(1-u)^{\lambda +\beta }}+\frac{1}{(1+u)^{\lambda
+\beta }}]u^{\beta +\sigma +\frac{2\varepsilon }{q}-1}du  \notag \\
&&+\int_{1}^{\infty }[\frac{1}{(u-1)^{\lambda +\beta }}+\frac{1}{%
(1+u)^{\lambda +\beta }}]u^{\sigma -\frac{2\varepsilon }{p}-1}du=\varepsilon
\tilde{I}<\varepsilon k\tilde{L}=k.  \label{24}
\end{eqnarray}

By (\ref{8}), (\ref{9}) and Fatou's lemma (cf. \cite{JK2}), we have
\begin{eqnarray*}
K(\sigma ) &=&\int_{0}^{1}[\frac{1}{(1-u)^{\lambda +\beta }}+\frac{1}{%
(1+u)^{\lambda +\beta }}]u^{\beta +\sigma -1}du \\
&&+\int_{1}^{\infty }[\frac{1}{(u-1)^{\lambda +\beta }}+\frac{1}{%
(1+u)^{\lambda +\beta }}]u^{\sigma -1}du
\end{eqnarray*}
\begin{eqnarray*}
&=&\int_{0}^{1}\lim_{\varepsilon \rightarrow 0^{+}}[\frac{1}{(1-u)^{\lambda
+\beta }}+\frac{1}{(1+u)^{\lambda +\beta }}]u^{\beta +\sigma +\frac{%
2\varepsilon }{q}-1}du \\
&&+\int_{1}^{\infty }\lim_{\varepsilon \rightarrow 0^{+}}[\frac{1}{%
(u-1)^{\lambda +\beta }}+\frac{1}{(1+u)^{\lambda +\beta }}]u^{\sigma -\frac{%
2\varepsilon }{p}-1}du \\
&\leq &\underset{\varepsilon \rightarrow 0^{+}}{\underline{\lim }}\left\{
\int_{0}^{1}[\frac{1}{(1-u)^{\lambda +\beta }}+\frac{1}{(1+u)^{\lambda
+\beta }}]u^{\beta +\sigma +\frac{2\varepsilon }{q}-1}du\right.
\end{eqnarray*}%
\begin{equation*}
\left. +\int_{1}^{\infty }[\frac{1}{(u-1)^{\lambda +\beta }}+\frac{1}{%
(1+u)^{\lambda +\beta }}]u^{\sigma -\frac{2\varepsilon }{p}-1}du\right\}
\leq k,
\end{equation*}%
which contradicts the fact that $k<K(\sigma )$. Hence the constant factor $%
K(\sigma )$ in (\ref{16}) is the best possible.

If the constant factor in (\ref{17}) is not the best possible, then by (\ref%
{20}) we would reach the contradiction that the constant factor in (\ref{16})
is not the best possible.
\end{proof}
\begin{theorem}
If in the assumptions of Theorem 3.1, we replace $p>1$ by $%
0<p<1$, we obtain the equivalent reverses of (\ref{16}) and (\ref{17}) with
the same best constant factors.
\end{theorem}
\begin{proof}
 By H\"{o}lder's reverse inequality (cf.
\cite{JK1}), we derive the reverses of (\ref{14}), (\ref{15}), (\ref{13}) and (%
\ref{20}). It is easy to obtain the reverse of (\ref{17}). In view of the
reverses of (\ref{17}) and (\ref{20}), we obtain the reverse of (\ref{16}).
On the other hand, if we suppose that the reverse of (\ref{16}) is valid, then if we set
$g(y)$ as in (\ref{21}), by the reverse of (\ref{15}), we have $J>0$.
If $J=\infty $, then the reverse of (\ref{17}) is obviously true; if $%
J<\infty $, then by the reverse of (\ref{16}), we obtain the reverses of (%
\ref{22}) and (\ref{23}). Hence, we obtain the reverse of (\ref{17}), which
is equivalent to the reverse of (\ref{16}).

If the constant factor $K(\sigma )$ in the reverse of (\ref{16}) is not the
best possible, then, there exists a positive constant $k,$ with $k>K(\sigma
) $, such that the reverse of (\ref{16}) is still valid when replacing $%
K(\sigma )$ by $k$. By the reverse of (\ref{24}), we have
\begin{equation*}
\int_{0}^{1}[\frac{1}{(1-u)^{\lambda +\beta }}+\frac{1}{(1+u)^{\lambda
+\beta }}]u^{\beta +\sigma +\frac{2\varepsilon }{q}-1}du
\end{equation*}%
\begin{equation}
+\int_{1}^{\infty }[\frac{1}{(u-1)^{\lambda +\beta }}+\frac{1}{%
(1+u)^{\lambda +\beta }}]u^{\sigma -\frac{2\varepsilon }{p}-1}du>k.
\label{25}
\end{equation}

By Levi's theorem (cf. \cite{JK2}), we find
\begin{eqnarray*}
&&\int_{1}^{\infty }[\frac{1}{(u-1)^{\lambda +\beta }}+\frac{1}{%
(1+u)^{\lambda +\beta }}]u^{\sigma -\frac{2\varepsilon }{p}-1}du \\
&\rightarrow &\int_{1}^{\infty }[\frac{1}{(u-1)^{\lambda +\beta }}+\frac{1}{%
(1+u)^{\lambda +\beta }}]u^{\sigma -1}du^{{}}(\varepsilon \rightarrow 0^{+}).
\end{eqnarray*}

There exists a constant $\delta _{0}>0,$ such that $\sigma -\frac{1}{2}%
\delta _{0}>-\beta ,$ and then $K(\sigma -\frac{\delta _{0}}{2})\in \mathbf{R%
}_{+}$. For $0<\varepsilon <\frac{\delta _{0}|q|}{4}\: (q<0)$, since $u^{\beta
+\sigma +\frac{2\varepsilon }{q}-1}\leq u^{\beta +\sigma -\frac{\delta _{0}}{%
2}-1},u\in (0,1],$ and
\begin{equation*}
0<\int_{0}^{1}[\frac{1}{(1-u)^{\lambda +\beta }}+\frac{1}{(1+u)^{\lambda
+\beta }}]u^{\beta +\sigma -\frac{\delta _{0}}{2}-1}du\leq K(\sigma -\frac{%
\delta _{0}}{2}),
\end{equation*}%
then by Lebesgue's control convergence theorem (cf. \cite{JK2}), we have
\begin{eqnarray*}
&&\int_{0}^{1}[\frac{1}{(1-u)^{\lambda +\beta }}+\frac{1}{(1+u)^{\lambda
+\beta }}]u^{\beta +\sigma +\frac{2\varepsilon }{q}-1}du \\
&\rightarrow &\int_{0}^{1}[\frac{1}{(1-u)^{\lambda +\beta }}+\frac{1}{%
(1+u)^{\lambda +\beta }}]u^{\beta +\sigma -1}du^{{}}\text{ (}\varepsilon
\rightarrow 0^{+}\text{)}.
\end{eqnarray*}%
By (\ref{25}) and the above results, for $\varepsilon \rightarrow 0^{+}$, we
get $K(\sigma )\geq k$, which contradicts the fact that $k>K(\sigma )$.
Hence, the constant factor $K(\sigma )$ in the reverse of (\ref{16}) is the
best possible.

If the constant factor in the reverse of (\ref{17}) is not the best possible,
then, by the reverse of (\ref{20}), we would reach the contradiction that the
constant factor in the reverse of (\ref{16}) is not the best possible.
\end{proof}
\begin{remark}
 (i) For $\delta =-1$ in (\ref{16}) and (\ref{17}),
replacing $|x|^{\lambda }f(x)$ by $f(x)$, we obtain $0<\int_{-\infty
}^{\infty }|x|^{p(1-\mu )-1}f^{p}(x)dx<\infty ,$ and the following
equivalent inequalities with the homogeneous kernel and the best possible
constant factors:
\begin{eqnarray}
&&\int_{-\infty }^{\infty }\int_{-\infty }^{\infty }\frac{(\min
\{|x|,|y|\})^{\beta }}{|x+y|^{\lambda +\beta }}f(x)g(y)dxdy  \notag \\
&<&K(\sigma )\left[ \int_{-\infty }^{\infty }|x|^{p(1-\mu )-1}f^{p}(x)dx%
\right] ^{\frac{1}{p}}\left[ \int_{-\infty }^{\infty }|y|^{q(1-\sigma
)-1}g^{q}(y)dy\right] ^{\frac{1}{q}},  \label{26}
\end{eqnarray}%
\begin{eqnarray}
&&\int_{-\infty }^{\infty }|y|^{p\sigma -1}\left[ \int_{-\infty }^{\infty }%
\frac{(\max \{|x|,|y|\})^{\beta }}{|x+y|^{\lambda +\beta }}f(x)dx\right]
^{p}dy  \notag \\
&<&K^{p}(\sigma )\int_{-\infty }^{\infty }|x|^{p(1-\mu )-1}f^{p}(x)dx.
\label{27}
\end{eqnarray}%
In particular, for $\lambda =0=\mu +\sigma (\mu ,\sigma >-\beta ),\:0<\beta <1,
$ we find
\begin{eqnarray}
K(\sigma ) &=&K_{0}(\sigma ):=\frac{1}{\beta +\sigma }F(\beta ,\beta +\sigma
,1+\beta +\sigma ,-1)  \notag \\
&&+\frac{1}{\beta +\mu }F(\beta ,\beta +\mu ,1+\beta +\mu ,-1)  \notag \\
&&+B(1-\beta ,\beta +\sigma )+B(1-\beta ,\beta +\mu ),  \label{28}
\end{eqnarray}%
and the following equivalent inequalities:%
\begin{eqnarray}
&&\int_{-\infty }^{\infty }\int_{-\infty }^{\infty }\left( \frac{\min
\{|x|,|y|\}}{|x+y|}\right) ^{\beta }f(x)g(y)dxdy  \notag \\
&<&K_{0}(\sigma )\left[ \int_{-\infty }^{\infty }|x|^{p(1-\mu )-1}f^{p}(x)dx%
\right] ^{\frac{1}{p}}\left[ \int_{-\infty }^{\infty }|y|^{q(1-\sigma
)-1}g^{q}(y)dy\right] ^{\frac{1}{q}},  \label{29}
\end{eqnarray}%
\begin{eqnarray}
&&\hspace{-3cm}\int_{-\infty }^{\infty }|y|^{p\sigma -1}\left[ \int_{-\infty
}^{\infty }\left( \frac{\min \{|x|,|y|\}}{|x+y|}\right) ^{\beta }f(x)dx%
\right] ^{p}dy  \notag \\
&<&K_{0}^{p}(\sigma )\int_{-\infty }^{\infty }|x|^{p(1-\mu )-1}f^{p}(x)dx.
\label{30}
\end{eqnarray}

(ii) For $\lambda =0=\mu +\sigma (\mu ,\sigma >-\beta ),\:0<\beta <1$ in (\ref%
{16}) and (\ref{17}), we have the following equivalent inequalities:%
\begin{eqnarray}
&&\int_{-\infty }^{\infty }\int_{-\infty }^{\infty }\left( \frac{\min
\{1,|x^{\delta }y|\}}{|1+x^{\delta }y|}\right) ^{\beta }f(x)g(y)dxdy  \notag
\\
&<&K_{0}(\sigma )\left[ \int_{-\infty }^{\infty }|x|^{p(1-\delta \sigma
)-1}f^{p}(x)dx\right] ^{\frac{1}{p}}\left[ \int_{-\infty }^{\infty
}|y|^{q(1-\sigma )-1}g^{q}(y)dy\right] ^{\frac{1}{q}},  \label{31}
\end{eqnarray}%
\begin{eqnarray}
&&\int_{-\infty }^{\infty }|y|^{p\sigma -1}\left[ \int_{-\infty }^{\infty
}\left( \frac{\min \{1,|x^{\delta }y|\}}{|1+x^{\delta }y|}\right) ^{\beta
}f(x)dx\right] ^{p}dy  \notag \\
&<&K_{0}^{p}(\sigma )\int_{-\infty }^{\infty }|x|^{p(1-\delta \sigma
)-1}f^{p}(x)dx.  \label{32}
\end{eqnarray}%
In particular, for $\delta =1,$we have the following equivalent inequalities
(cf. \cite{BHY2}, for \mbox{$\sigma =\mu =0$}):
\begin{eqnarray}
&&\int_{-\infty }^{\infty }\int_{-\infty }^{\infty }\left( \frac{\min
\{1,|xy|\}}{|1+xy|}\right) ^{\beta }f(x)g(y)dxdy  \notag \\
&<&K_{0}(\sigma )\left[ \int_{-\infty }^{\infty }|x|^{p(1-\sigma
)-1}f^{p}(x)dx\right] ^{\frac{1}{p}}\left[ \int_{-\infty }^{\infty
}|y|^{q(1-\sigma )-1}g^{q}(y)dy\right] ^{\frac{1}{q}},  \label{33}
\end{eqnarray}%
\begin{eqnarray}
&&\int_{-\infty }^{\infty }|y|^{p\sigma -1}\left[ \int_{-\infty }^{\infty
}\left( \frac{\min \{1,|xy|\}}{|1+xy|}\right) ^{\beta }f(x)dx\right] ^{p}dy
\notag \\
&<&K_{0}^{p}(\sigma )\int_{-\infty }^{\infty }|x|^{p(1-\sigma )-1}f^{p}(x)dx.
\label{34}
\end{eqnarray}

(iii) For $\beta =0<\lambda <1,\:\sigma =\mu =\frac{\lambda }{2}$ in (\ref{18}%
), we obtain%
\begin{equation*}
K(\frac{\lambda }{2})=\int_{0}^{\infty }\frac{u^{\frac{\lambda }{2}-1}}{%
(1+u)^{\lambda }}du+2\int_{0}^{1}\frac{u^{\frac{\lambda }{2}-1}}{%
(1-u)^{\lambda }}du=k_{\lambda },
\end{equation*}%
and then (\ref{3}) follows. Hence, (\ref{16})-(\ref{18}) is an extension
of (\ref{3}).
\end{remark}
\section{Some Corollaries}
In the following two sections, if the constant factors are related to $%
K_{1}(\sigma ),$ then we call them Hardy-type inequalities (operators) of the first kind; if the constant factors are related to $K_{2}(\sigma ),$ then we
call them Hardy-type inequalities (operators) of the second kind.

Setting the kernel%
\begin{equation*}
H(xy):=\left\{
\begin{array}{c}
0,\ \ \ \ \ \ \ \ \ \ \ \ \ \ \ \ \ \ |xy|>1 \\
\frac{(\min \{1,|xy|\})^{\beta }}{|1+xy|^{\lambda +\beta }},\ \ |xy|\leq 1%
\end{array}%
\right. ,
\end{equation*}%
it follows that
\begin{eqnarray*}
H(u) &=&\left\{
\begin{array}{c}
0,\ \ \ \ \ \ \ \ \ \ \ \ \ \ \ \ \ \ |u|>1 \\
\frac{(\min \{1,|u|\})^{\beta }}{|1+u|^{\lambda +\beta }},\ \ \ |u|\leq 1%
\end{array}%
\right. , \\
\int_{-\infty }^{\infty }H(u)|u|^{\sigma -1}du &=&\int_{-1}^{1}\frac{(\min
\{1,|u|\})^{\beta }}{|1+u|^{\lambda +\beta }}|u|^{\sigma -1}du=K_{1}(\sigma
).
\end{eqnarray*}%
In view of Theorems 3.1-3.2 (for $\delta =1)$, we obtain the following 
Hardy-type inequalities of the first kind with the non-homogeneous kernel:
\begin{corollary}
Suppose that $p>1,\frac{1}{p}+\frac{1}{q}=1,\beta
>-1,\mu ,\sigma >-\beta ,$ $\mu +\sigma =\lambda <1-\beta ,$ $K_{1}(\sigma )$
is indicated by (\ref{8}) (or \ref{10}). If $f,\:g\geq 0,$ satisfy $$%
0<\int_{-\infty }^{\infty }|x|^{p(1-\sigma )-1}f^{p}(x)dx<\infty $$ and $$%
0<\int_{-\infty }^{\infty }|y|^{q(1-\sigma )-1}g^{q}(y)dy<\infty ,$$ then we
have the following equivalent inequalities:
\begin{eqnarray}
&&\int_{-\infty }^{\infty }\left[ \int_{-\frac{1}{|y|}}^{\frac{1}{|y|}}\frac{%
(\min \{1,|xy|\})^{\beta }}{|1+xy|^{\lambda +\beta }}f(x)dx\right] g(y)dy
\notag \\
&<&K_{1}(\sigma )\left[ \int_{-\infty }^{\infty }|x|^{p(1-\sigma
)-1}f^{p}(x)dx\right] ^{\frac{1}{p}}\left[ \int_{-\infty }^{\infty
}|y|^{q(1-\sigma )-1}g^{q}(y)dy\right] ^{\frac{1}{q}},  \label{35}
\end{eqnarray}%
\begin{eqnarray}
&&\int_{-\infty }^{\infty }|y|^{p\sigma -1}\left[ \int_{-\frac{1}{|y|}}^{%
\frac{1}{|y|}}\frac{(\min \{1,|xy|\})^{\beta }}{|1+x^{\delta }y|^{\lambda
+\beta }}f(x)dx\right] ^{p}dy  \notag \\
&<&K_{1}^{p}(\sigma )\int_{-\infty }^{\infty }|x|^{p(1-\sigma )-1}f^{p}(x)dx.
\label{36}
\end{eqnarray}%
where, the constant factors $K_{1}(\sigma )$ and $K_{1}^{p}(\sigma )$ are the
best possible. Replacing $p>1$ by $0<p<1$, we have the equivalent reverses
of (\ref{35}) and (\ref{36}) with the same best constant factors.
\end{corollary}

If we set $E_{y}:=\{x\in \mathbf{R};|xy|\geq 1\},$ and%
\begin{equation*}
H(xy):=\left\{
\begin{array}{c}
0,\ \ \ \ \ \ \ \ \ \ \ \ \ \ \ \ \ \ |xy|<1 \\
\frac{(\min \{1,|xy|\})^{\beta }}{|1+xy|^{\lambda +\beta }},\ \ |xy|\geq 1%
\end{array}%
\right. ,
\end{equation*}%
then it follows that
\begin{eqnarray*}
H(u) &=&\left\{
\begin{array}{c}
0,\ \ \ \ \ \ \ \ \ \ \ \ \ \ \ \ \ \ |u|<1 \\
\frac{(\min \{1,|u|\})^{\beta }}{|1+u|^{\lambda +\beta }},\ \ \ |u|\geq 1%
\end{array}%
\right. , \\
\int_{-\infty }^{\infty }H(u)|u|^{\sigma -1}du &=&\int_{E_{1}}\frac{(\min
\{1,|u|\})^{\beta }}{|1+u|^{\lambda +\beta }}|u|^{\sigma -1}du=K_{2}(\sigma
).
\end{eqnarray*}%
In view of Theorems 3.1-3.2 (for $\delta =1)$, we have the following 
Hardy-type inequalities of the second kind with the non-homogeneous kernel:
\begin{corollary}
Suppose that $p>1,\frac{1}{p}+\frac{1}{q}=1,\beta
>-1,\mu ,\sigma >-\beta ,$ $\mu +\sigma =\lambda <1-\beta ,$ $K_{2}(\sigma )$
is indicated by (\ref{9}) (or \ref{11}).\\ 
If $f,\:g\geq 0,$ satisfy $$%
0<\int_{-\infty }^{\infty }|x|^{p(1-\sigma )-1}f^{p}(x)dx<\infty $$ and $$%
0<\int_{-\infty }^{\infty }|y|^{q(1-\sigma )-1}g^{q}(y)dy<\infty, $$ then we
have the following equivalent inequalities:
\begin{eqnarray}
&&\int_{-\infty }^{\infty }\left[ \int_{E_{y}}\frac{(\min \{1,|xy|\})^{\beta
}}{|1+xy|^{\lambda +\beta }}f(x)dx\right] g(y)dy  \notag \\
&<&K_{2}(\sigma )\left[ \int_{-\infty }^{\infty }|x|^{p(1-\sigma
)-1}f^{p}(x)dx\right] ^{\frac{1}{p}}\left[ \int_{-\infty }^{\infty
}|y|^{q(1-\sigma )-1}g^{q}(y)dy\right] ^{\frac{1}{q}},  \label{37}
\end{eqnarray}%
\begin{eqnarray}
&&\int_{-\infty }^{\infty }|y|^{p\sigma -1}\left[ \int_{E_{y}}\frac{(\min
\{1,|xy|\})^{\beta }}{|1+xy|^{\lambda +\beta }}f(x)dx\right] ^{p}dy  \notag
\\
&<&K_{2}^{p}(\sigma )\int_{-\infty }^{\infty }|x|^{p(1-\delta \sigma
)-1}f^{p}(x)dx.  \label{38}
\end{eqnarray}%
where, the constant factors $K_{2}(\sigma )$ and $K_{2}^{p}(\sigma )$ are the
best possible. Replacing $p>1$ by $0<p<1$, we obtain the equivalent reverses
of (\ref{37}) and (\ref{38}) with the same best constant factors.
\end{corollary}

If we set $\widetilde{E}_{y}:=\{x\in \mathbf{R};|\frac{y}{x}|\leq 1\}$ and%
\begin{equation*}
K_{\lambda }(x,y):=\left\{
\begin{array}{c}
0,\ \ \ \ \ \ \ \ \ \ \ \ \ \ \ \ \ \ |\frac{y}{x}|>1 \\
\frac{(\min \{|x|,|y|\})^{\beta }}{|x+y|^{\lambda +\beta }},\ \ |\frac{y}{x}%
|\leq 1%
\end{array}%
\right. ,
\end{equation*}%
then it follows
\begin{eqnarray*}
K_{\lambda }(1,u) &=&\left\{
\begin{array}{c}
0,\ \ \ \ \ \ \ \ \ \ \ \ \ \ \ \ \ \ |u|>1 \\
\frac{(\min \{1,|u|\})^{\beta }}{|1+u|^{\lambda +\beta }},\ \ \ |u|\leq 1%
\end{array}%
\right. , \\
\int_{-\infty }^{\infty }K_{\lambda }(1,u)|u|^{\sigma -1}du &=&\int_{-1}^{1}%
\frac{(\min \{1,|u|\})^{\beta }}{|1+u|^{\lambda +\beta }}|u|^{\sigma
-1}du=K_{1}(\sigma ).
\end{eqnarray*}%
In view of Remark 3.3 (i), we have the following Hardy-type
inequalities of the first kind with the homogeneous kernel:
\begin{corollary}
Suppose that $p>1,\frac{1}{p}+\frac{1}{q}=1,\beta
>-1,\mu ,\sigma >-\beta ,$ $\mu +\sigma =\lambda <1-\beta ,$ $K_{1}(\sigma )$
is indicated by (\ref{8}) (or \ref{10}). If $f,\:g\geq 0,$ satisfy $$%
0<\int_{-\infty }^{\infty }|x|^{p(1-\mu )-1}f^{p}(x)dx<\infty $$ and $$%
0<\int_{-\infty }^{\infty }|y|^{q(1-\sigma )-1}g^{q}(y)dy<\infty ,$$ then we
have the following equivalent inequalities:
\begin{eqnarray}
&&\int_{-\infty }^{\infty }\left[ \int_{\widetilde{E}_{y}}\frac{(\min
\{|x|,|y|\})^{\beta }}{|x+y|^{\lambda +\beta }}f(x)dx\right] g(y)dy  \notag
\\
&<&K_{1}(\sigma )\left[ \int_{-\infty }^{\infty }|x|^{p(1-\mu )-1}f^{p}(x)dx%
\right] ^{\frac{1}{p}}\left[ \int_{-\infty }^{\infty }|y|^{q(1-\sigma
)-1}g^{q}(y)dy\right] ^{\frac{1}{q}},  \label{39}
\end{eqnarray}%
\begin{eqnarray}
&&\int_{-\infty }^{\infty }|y|^{p\sigma -1}\left[ \int_{\widetilde{E}_{y}}%
\frac{(\min \{|x|,|y|\})^{\beta }}{|x+y|^{\lambda +\beta }}f(x)dx\right]
^{p}dy  \notag \\
&<&K_{1}^{p}(\sigma )\int_{-\infty }^{\infty }|x|^{p(1-\mu )-1}f^{p}(x)dx.
\label{40}
\end{eqnarray}%
where the constant factors $K_{1}(\sigma )$ and $K_{1}^{p}(\sigma )$ are the
best possible. Replacing $p>1$ by $0<p<1$, we derive the equivalent reverses
of (\ref{39}) and (\ref{40}) with the same best constant factors.
\end{corollary}

Setting the kernel
\begin{equation*}
K_{\lambda }(x,y):=\left\{
\begin{array}{c}
0,\ \ \ \ \ \ \ \ \ \ \ \ \ \ \ \ \ \ |\frac{y}{x}|<1 \\
\frac{(\min \{|x|,|y|\})^{\beta }}{|x+y|^{\lambda +\beta }},\ \ |\frac{y}{x}%
|\geq 1%
\end{array}%
\right. ,
\end{equation*}%
then it follows that
\begin{eqnarray*}
K_{\lambda }(1,u) &=&\left\{
\begin{array}{c}
0,\ \ \ \ \ \ \ \ \ \ \ \ \ \ \ \ \ \ |u|<1 \\
\frac{(\min \{1,|u|\})^{\beta }}{|1+u|^{\lambda +\beta }},\ \ \ |u|\geq 1%
\end{array}%
\right. , \\
\int_{-\infty }^{\infty }K_{\lambda }(1,u)|u|^{\sigma -1}du &=&\int_{E_{1}}%
\frac{(\min \{1,|u|\})^{\beta }}{|1+u|^{\lambda +\beta }}|u|^{\sigma
-1}du=K_{2}(\sigma ).
\end{eqnarray*}%
In view of Remark 3.3 (i), we have the following Hardy-type
inequalities of the second kind with the homogeneous kernel:
\begin{corollary}
Suppose that $p>1,\frac{1}{p}+\frac{1}{q}=1,\beta
<1,\mu ,\sigma >0,$ $\mu +\sigma =\lambda <1-\beta ,K_{2}(\sigma )$ is
indicated by (\ref{9}) (or \ref{11}). If $f,\:g\geq 0,$ satisfy $$%
0<\int_{-\infty }^{\infty }|x|^{p(1-\mu )-1}f^{p}(x)dx<\infty $$ and $$%
0<\int_{-\infty }^{\infty }|y|^{q(1-\sigma )-1}g^{q}(y)dy<\infty ,$$ then we
have the following equivalent inequalities:
\begin{eqnarray}
&&\int_{-\infty }^{\infty }\left[ \int_{-|y|}^{|y|}\frac{(\min
\{|x|,|y|\})^{\beta }}{|x+y|^{\lambda +\beta }}f(x)dx\right] g(y)dy  \notag
\\
&<&K_{2}(\sigma )\left[ \int_{-\infty }^{\infty }|x|^{p(1-\mu )-1}f^{p}(x)dx%
\right] ^{\frac{1}{p}}\left[ \int_{-\infty }^{\infty }|y|^{q(1-\sigma
)-1}g^{q}(y)dy\right] ^{\frac{1}{q}},  \label{41}
\end{eqnarray}%
\begin{eqnarray}
&&\int_{-\infty }^{\infty }|y|^{p\sigma -1}\left[ \int_{-|y|}^{|y|}\frac{%
(\min \{|x|,|y|\})^{\beta }}{|x+y|^{\lambda +\beta }}f(x)dx\right] ^{p}dy
\notag \\
&<&K_{2}^{p}(\sigma )\int_{-\infty }^{\infty }|x|^{p(1-\mu )-1}f^{p}(x)dx.
\label{42}
\end{eqnarray}%
where the constant factors $K_{2}(\sigma )$ and $K_{2}^{p}(\sigma )$ are the
best possible. Replacing $p>1$ by $0<p<1$, we get the equivalent reverses
of (\ref{41}) and (\ref{42}) with the same best constant factors.
\end{corollary}
\section{Operator Expressions}
Suppose that $p>1,\frac{1}{p}+\frac{1}{q}=1,\beta >-1,\mu ,\sigma >-\beta ,$
$\mu +\sigma =\lambda <1-\beta .$ We set the following functions: $\varphi
(x):=|x|^{p(1-\sigma )-1},\psi (y):=|y|^{q(1-\sigma )-1},\phi
(x):=|x|^{p(1-\mu )-1}(x,y\in \mathbf{R}),$ wherefrom, $\psi
^{1-p}(y)=|y|^{p\sigma -1}.$ Define the following real normed linear space:%
\begin{equation*}
L_{p,\varphi }(\mathbf{R}):=\left\{ f:||f||_{p,\varphi }:=\left(
\int_{-\infty }^{\infty }\varphi (x)|f(x)|^{p}dx\right) ^{\frac{1}{p}%
}<\infty \right\} ,
\end{equation*}%
wherefrom,%
\begin{eqnarray*}
L_{p,\psi ^{1-p}}(\mathbf{R}) &=&\left\{ h:||h||_{p,\psi ^{1-p}}=\left(
\int_{-\infty }^{\infty }\psi ^{1-p}(y)|h(y)|^{p}dy\right) ^{\frac{1}{p}%
}<\infty \right\} , \\
L_{p,\phi }(\mathbf{R}) &=&\left\{ g:||g||_{p,\phi }=\left( \int_{-\infty
}^{\infty }\phi (x)|g(x)|^{p}dx\right) ^{\frac{1}{p}}<\infty \right\} .
\end{eqnarray*}

(a) In view of Theorem 3.1 ($\delta =1)$, for $f\in L_{p,\varphi }(\mathbf{R}),$
setting%
\begin{equation*}
H^{(1)}(y):=\int_{-\infty }^{\infty }\frac{(\min \{1,|xy|\})^{\beta }}{%
|1+xy|^{\lambda +\beta }}|f(x)|dx^{{}}(y\in \mathbf{R}),
\end{equation*}%
by (\ref{19}), we have
\begin{equation}
||H^{(1)}||_{p,\psi ^{1-p}}=\left[ \int_{-\infty }^{\infty }\psi
^{1-p}(y)(H^{(1)}(y))^{p}dy\right] ^{\frac{1}{p}}<K(\sigma )||f||_{p,\varphi
}<\infty .  \label{43}
\end{equation}
\begin{definition}
 Define the Hilbert-type integral operator with the
non-homogeneous kernel in the whole plane $T^{(1)}\::\:L_{p,\varphi }(\mathbf{%
R})\rightarrow L_{p,\psi ^{1-p}}(\mathbf{R})$ as follows: For any $f\in
L_{p,\varphi }(\mathbf{R}),$ there exists a unique representation $%
T^{(1)}f=H^{(1)}\in L_{p,\psi ^{1-p}}(\mathbf{R}),$ satisfying $$T^{(1)}f(y)=H^{(1)}(y),$$
for any $y\in
\mathbf{R}$.
\end{definition}

In view of (\ref{43}), it follows that $||T^{(1)}f||_{p,\psi
^{1-p}}=||H^{(1)}||_{p,\psi ^{1-p}}\leq K(\sigma )||f||_{p,\varphi }.$ 
Then, the operator $T^{(1)}$ is bounded satisfying%
\begin{equation*}
||T^{(1)}||=\sup_{f(\neq \theta )\in L_{p,\varphi }(\mathbf{R}_{+})}\frac{%
||T^{(1)}f||_{p,\psi ^{1-p}}}{||f||_{p,\varphi }}\leq K(\sigma ).
\end{equation*}%
Since the constant factor $K(\sigma )$ in (\ref{43}) is the best possible,
we have $||T^{(1)}||=K(\sigma ).$

If we define the formal inner product of $T^{(1)}f$ and $g$ as follows:
\begin{eqnarray*}
(T^{(1)}f,g) &:&=\int_{-\infty }^{\infty }\left[ \int_{-\infty }^{\infty }%
\frac{(\min \{1,|xy|\})^{\beta }}{|1+xy|^{\lambda +\beta }}f(x)dx\right]
g(y)dy \\
&=&\int_{-\infty }^{\infty }\int_{-\infty }^{\infty }\frac{(\min
\{1,|xy|\})^{\beta }}{|1+xy|^{\lambda +\beta }}f(x)g(y)dxdy,
\end{eqnarray*}%
then we can rewrite (\ref{18}) and (\ref{19}) in the form:%
\begin{equation}
(T^{(1)}f,g)<||T^{(1)}||\cdot ||f||_{p,\varphi }||g||_{q,\psi
},\ ||T^{(1)}f||_{p,\psi ^{1-p}}<||T^{(1)}||\cdot ||f||_{p,\varphi }.
\label{44}
\end{equation}

(b) In view of Corollary 4.1, for $f\in L_{p,\varphi }(\mathbf{R}),$ setting%
\begin{equation*}
H_{1}^{(1)}(y):=\int_{-\frac{1}{|y|}}^{\frac{1}{|y|}}\frac{(\min
\{1,|xy|\})^{\beta }}{|1+xy|^{\lambda +\beta }}|f(x)|dx^{{}}(y\in \mathbf{R}%
),
\end{equation*}%
by (\ref{29}), we obtain
\begin{equation}
||H_{1}^{(1)}||_{p,\psi ^{1-p}}=\left[ \int_{-\frac{1}{|y|}}^{\frac{1}{|y|}%
}\psi ^{1-p}(y)(H_{1}^{(1)}(y))^{p}dy\right] ^{\frac{1}{p}}<K_{1}(\sigma
)||f||_{p,\varphi }<\infty .  \label{45}
\end{equation}
\begin{definition}
 Define the Hilbert-type integral operator of the first kind
with the non-homogeneous kernel in the whole plane $T_{1}^{(1)}:$ $%
L_{p,\varphi }(\mathbf{R})\rightarrow L_{p,\psi ^{1-p}}(\mathbf{R})$ as
follows: For any $f\in L_{p,\varphi }(\mathbf{R}),$ there exists a unique
representation $T_{1}^{(1)}f=H_{1}^{(1)}\in L_{p,\psi ^{1-p}}(\mathbf{R}),$
satisfying $$T_{1}^{(1)}f(y)=H_{1}^{(1)}(y),$$
for any $y\in \mathbf{R}$.
\end{definition}

In view of (\ref{45}), it follows that $||T_{1}^{(1)}f||_{p,\psi
^{1-p}}=||H_{1}^{(1)}||_{p,\psi ^{1-p}}\leq K_{1}(\sigma )||f||_{p,\varphi
}, $ and then the operator $T_{1}^{(1)}$ is bounded satisfying%
\begin{equation*}
||T_{1}^{(1)}||=\sup_{f(\neq \theta )\in L_{p,\varphi }(\mathbf{R}_{+})}%
\frac{||T_{1}^{(1)}f||_{p,\psi ^{1-p}}}{||f||_{p,\varphi }}\leq K_{1}(\sigma
).
\end{equation*}%
Since the constant factor $K_{1}(\sigma )$ in (\ref{45}) is the best
possible, we have $||T_{1}^{(1)}||=K_{1}(\sigma ).$

If we define the formal inner product of $T_{1}^{(1)}f$ and $g$ as follows:
\begin{equation*}
(T_{1}^{(1)}f,g):=\int_{-\infty }^{\infty }\left[ \int_{-\frac{1}{|y|}}^{%
\frac{1}{|y|}}\frac{(\min \{1,|xy|\})^{\beta }}{|1+xy|^{\lambda +\beta }}%
f(x)dx\right] g(y)dy,
\end{equation*}%
then we can rewrite (\ref{28}) and (\ref{29}) in the following way:%
\begin{equation}
(T_{1}^{(1)}f,g)<||T_{1}^{(1)}||\cdot ||f||_{p,\varphi }||g||_{q,\psi
},||T_{1}^{(1)}f||_{p,\psi ^{1-p}}<||T_{1}^{(1)}||\cdot ||f||_{p,\varphi }.
\label{46}
\end{equation}

(c) In view of Corollary 4.2, for $f\in L_{p,\varphi }(\mathbf{R}),$ setting%
\begin{equation*}
H_{2}^{(1)}(y):=\int_{E_{y}}\frac{(\min \{1,|xy|\})^{\beta }}{%
|1+xy|^{\lambda +\beta }}|f(x)|dx^{{}}(y\in \mathbf{R}),
\end{equation*}%
by (\ref{31}), we have
\begin{equation}
||H_{2}^{(1)}||_{p,\psi ^{1-p}}=\left[ \int_{E_{y}}\psi
^{1-p}(y)(H_{2}^{(1)}(y))^{p}dy\right] ^{\frac{1}{p}}<K_{2}(\sigma
)||f||_{p,\varphi }<\infty .  \label{47}
\end{equation}
\begin{definition}
 Define the Hilbert-type integral operator of the second kind
with the non-homogeneous kernel in the whole plane $T_{2}^{(1)}\::\:
L_{p,\varphi }(\mathbf{R})\rightarrow L_{p,\psi ^{1-p}}(\mathbf{R})$ as
follows: For any $f\in L_{p,\varphi }(\mathbf{R}),$ there exists a unique
representation $T_{2}^{(1)}f=H_{2}^{(1)}\in L_{p,\psi ^{1-p}}(\mathbf{R}),$
satisfying $$T_{2}^{(1)}f(y)=H_{2}^{(1)}(y),$$
for any $y\in \mathbf{R}.$
\end{definition}

In view of (\ref{47}), it follows that $||T_{2}^{(1)}f||_{p,\psi
^{1-p}}=||H_{2}^{(1)}||_{p,\psi ^{1-p}}\leq K_{2}(\sigma )||f||_{p,\varphi
}, $ and then the operator $T_{2}^{(1)}$ is bounded satisfying%
\begin{equation*}
||T_{2}^{(1)}||=\sup_{f(\neq \theta )\in L_{p,\varphi }(\mathbf{R}_{+})}%
\frac{||T_{2}^{(1)}f||_{p,\psi ^{1-p}}}{||f||_{p,\varphi }}\leq K_{2}(\sigma
).
\end{equation*}%
Since the constant factor $K_{2}(\sigma )$ in (\ref{47}) is the best
possible, we have $||T_{2}^{(1)}||=K_{2}(\sigma ).$

If we define the formal inner product of $T_{2}^{(1)}f$ and $g$ as follows:
\begin{equation*}
(T_{2}^{(1)}f,g):=\int_{-\infty }^{\infty }\left[ \int_{E_{y}}\frac{(\min
\{1,|xy|\})^{\beta }}{|1+xy|^{\lambda +\beta }}f(x)dx\right] g(y)dy,
\end{equation*}%
then we can rewrite (\ref{30}) and (\ref{31}) as shown below:%
\begin{equation}
(T_{2}^{(1)}f,g)<||T_{2}^{(1)}||\cdot ||f||_{p,\varphi }||g||_{q,\psi
},||T_{2}^{(1)}f||_{p,\psi ^{1-p}}<||T_{2}^{(1)}||\cdot ||f||_{p,\varphi }.
\label{48}
\end{equation}

(d) In view of Remark 3.3 (i), for $f\in L_{p,\phi }(\mathbf{R}),$ setting%
\begin{equation*}
H^{(2)}(y):=\int_{-\infty }^{\infty }\frac{(\min \{|x|,|y|\})^{\beta }}{%
|x+y|^{\lambda +\beta }}|f(x)|dx^{{}}(y\in \mathbf{R}),
\end{equation*}%
by (\ref{27}), we have
\begin{equation}
||H^{(2)}||_{p,\psi ^{1-p}}=\left[ \int_{-\infty }^{\infty }\psi
^{1-p}(y)(H^{(2)}(y))^{p}dy\right] ^{\frac{1}{p}}<K(\sigma )||f||_{p,\phi
}<\infty .  \label{49}
\end{equation}
\begin{definition}
Define the Hilbert-type integral operator with the
homogeneous kernel in the whole plane $T^{(2)}\::\:L_{p,\phi }(\mathbf{R}%
)\rightarrow L_{p,\psi ^{1-p}}(\mathbf{R})$ as follows: For any $f\in
L_{p,\phi }(\mathbf{R}),$ there exists a unique representation $%
T^{(2)}f=H^{(2)}\in L_{p,\psi ^{1-p}}(\mathbf{R}),$ satisfying $$T^{(2)}f(y)=H^{(2)}(y),$$
for any $y\in
\mathbf{R}.$
\end{definition}

In view of (\ref{49}), it follows that $||T^{(2)}f||_{p,\psi
^{1-p}}=||H^{(2)}||_{p,\psi ^{1-p}}\leq K(\sigma )||f||_{p,\phi },$ and then
the operator $T^{(2)}$ is bounded satisfying%
\begin{equation*}
||T^{(2)}||=\sup_{f(\neq \theta )\in L_{p,\phi }(\mathbf{R}_{+})}\frac{%
||T^{(2)}f||_{p,\psi ^{1-p}}}{||f||_{p,\phi }}\leq K(\sigma ).
\end{equation*}%
Since the constant factor $K(\sigma )$ in (\ref{49}) is the best possible,
we have $||T^{(2)}||=K(\sigma ).$

If we define the formal inner product of $T^{(2)}f$ and $g$ as follows:
\begin{eqnarray*}
(T^{(2)}f,g) &:&=\int_{-\infty }^{\infty }\left[ \int_{-\infty }^{\infty }%
\frac{(\min \{|x|,|y|\})^{\beta }}{|x+y|^{\lambda +\beta }}f(x)dx\right]
g(y)dy \\
&=&\int_{-\infty }^{\infty }\int_{-\infty }^{\infty }\frac{(\min
\{|x|,|y|\})^{\beta }}{|x+y|^{\lambda +\beta }}f(x)g(y)dxdy,
\end{eqnarray*}%
then we can rewrite (\ref{26}) and (\ref{27}) as follows:%
\begin{equation}
(T^{(2)}f,g)<||T^{(2)}||\cdot ||f||_{p,\phi }||g||_{q,\psi
},||T^{(2)}f||_{p,\psi ^{1-p}}<||T^{(2)}||\cdot ||f||_{p,\phi }.  \label{50}
\end{equation}

(e) In view of Corollary 4.3, for $f\in L_{p,\phi }(\mathbf{R}),$ setting%
\begin{equation*}
H_{1}^{(2)}(y):=\int_{\widetilde{E}_{y}}\frac{(\min \{|x|,|y|\})^{\beta }}{%
|x+y|^{\lambda +\beta }}|f(x)|dx^{{}}(y\in \mathbf{R}),
\end{equation*}%
by (\ref{33}), we have
\begin{equation}
||H_{1}^{(2)}||_{p,\psi ^{1-p}}=\left[ \int_{\widetilde{E}_{y}}\psi
^{1-p}(y)(H_{1}^{(2)}(y))^{p}dy\right] ^{\frac{1}{p}}<K_{1}(\sigma
)||f||_{p,\phi }<\infty .  \label{51}
\end{equation}
\begin{definition}
 Define the Hilbert-type integral operator of the fist kind
with the homogeneous kernel in the whole plane $T_{1}^{(2)}\::\:L_{p,\phi }(%
\mathbf{R})\rightarrow L_{p,\psi ^{1-p}}(\mathbf{R})$ as follows: For any $%
f\in L_{p,\phi }(\mathbf{R}),$ there exists a unique representation $%
T_{1}^{(2)}f=H_{1}^{(2)}\in L_{p,\psi ^{1-p}}(\mathbf{R}),$ satisfying 
$$T_{1}^{(2)}f(y)=H_{1}^{(2)}(y),$$
for
any $y\in \mathbf{R}.$
\end{definition}

In view of (\ref{51}), it follows that $||T_{1}^{(2)}f||_{p,\psi
^{1-p}}=||H_{1}^{(2)}||_{p,\psi ^{1-p}}\leq K_{1}(\sigma )||f||_{p,\phi },$
and then the operator $T_{1}^{(2)}$ is bounded satisfying%
\begin{equation*}
||T_{1}^{(2)}||=\sup_{f(\neq \theta )\in L_{p,\phi }(\mathbf{R}_{+})}\frac{%
||T_{1}^{(2)}f||_{p,\psi ^{1-p}}}{||f||_{p,\phi }}\leq K_{1}(\sigma ).
\end{equation*}%
Since the constant factor $K_{1}(\sigma )$ in (\ref{51}) is the best
possible, we have $||T_{1}^{(2)}||=K_{1}(\sigma ).$

If we define the formal inner product of $T_{1}^{(2)}f$ and $g$ as follows:
\begin{equation*}
(T_{1}^{(1)}f,g):=\int_{-\infty }^{\infty }\left[ \int_{\widetilde{E}_{y}}%
\frac{(\min \{|x|,|y|\})^{\beta }}{|x+y|^{\lambda +\beta }}f(x)dx\right]
g(y)dy,
\end{equation*}%
then we can rewrite (\ref{32}) and (\ref{33}) as follows:%
\begin{equation}
(T_{1}^{(2)}f,g)<||T_{1}^{(2)}||\cdot ||f||_{p,\phi }||g||_{q,\psi
},||T_{1}^{(2)}f||_{p,\psi ^{1-p}}<||T_{1}^{(2)}||\cdot ||f||_{p,\phi }.
\label{52}
\end{equation}

(f) In view of Corollary 4.4, for $f\in L_{p,\phi }(\mathbf{R}),$ setting%
\begin{equation*}
H_{2}^{(2)}(y):=\int_{-|y|}^{|y|}\frac{(\min \{|x|,|y|\})^{\beta }}{%
|x+y|^{\lambda +\beta }}|f(x)|dx^{{}}(y\in \mathbf{R}),
\end{equation*}%
by (\ref{35}), we have
\begin{equation}
||H_{2}^{(2)}||_{p,\psi ^{1-p}}=\left[ \int_{-|y|}^{|y|}\psi
^{1-p}(y)(H_{2}^{(2)}(y))^{p}dy\right] ^{\frac{1}{p}}<K_{2}(\sigma
)||f||_{p,\phi }<\infty .  \label{53}
\end{equation}
\begin{definition}
 Define the Hilbert-type integral operator of the second kind
with the homogeneous kernel in the whole plane $T_{2}^{(2)}\::\:L_{p,\phi }(%
\mathbf{R})\rightarrow L_{p,\psi ^{1-p}}(\mathbf{R})$ as follows: For any $%
f\in L_{p,\phi }(\mathbf{R}),$ there exists a unique representation $%
T_{2}^{(2)}f=H_{2}^{(2)}\in L_{p,\psi ^{1-p}}(\mathbf{R}),$ satisfying 
$$T_{2}^{(2)}f(y)=H_{2}^{(2)}(y),$$
for
any $y\in \mathbf{R}.$
\end{definition}

In view of (\ref{53}), it follows that $||T_{2}^{(2)}f||_{p,\psi
^{1-p}}=||H_{2}^{(2)}||_{p,\psi ^{1-p}}\leq K_{2}(\sigma )||f||_{p,\phi },$
and thus the operator $T_{2}^{(2)}$ is bounded satisfying%
\begin{equation*}
||T_{2}^{(2)}||=\sup_{f(\neq \theta )\in L_{p,\phi }(\mathbf{R}_{+})}\frac{%
||T_{2}^{(2)}f||_{p,\psi ^{1-p}}}{||f||_{p,\phi }}\leq K_{2}(\sigma ).
\end{equation*}%
Since the constant factor $K_{2}(\sigma )$ in (\ref{53}) is the best
possible, we have $||T_{2}^{(2)}||=K_{2}(\sigma ).$

If we define the formal inner product of $T_{2}^{(2)}f$ and $g$ as
\begin{equation*}
(T_{2}^{(2)}f,g):=\int_{-\infty }^{\infty }\left[ \int_{-|y|}^{|y|}\frac{%
(\min \{|x|,|y|\})^{\beta }}{|x+y|^{\lambda +\beta }}f(x)dx\right] g(y)dy,
\end{equation*}%
then we can rewrite (\ref{34}) and (\ref{35}) as follows:%
\begin{equation}
(T_{2}^{(2)}f,g)<||T_{2}^{(2)}||\cdot ||f||_{p,\phi }||g||_{q,\psi
},||T_{2}^{(2)}f||_{p,\psi ^{1-p}}<||T_{2}^{(2)}||\cdot ||f||_{p,\phi }.
\label{54}
\end{equation}

\section*{Acknowledgments}
The authors wish to express their thanks to Professors Mario Krnic and Jichang Kuang for their careful reading of the manuscript and their valuable suggestions.


\begin{thebibliography}{99}
\bibitem{HLP} G.H. Hardy, J.E. Littlewood, G.P$\acute{o}$lya, Inequalities,
Cambridge University Press, Cambridge, USA, 1934.

\bibitem{MPF} D.S. Mitrinovi$\acute{c}$, J.E. Pecaric, A.M.Fink,
Inequalities Involving Functions and their Integrals and Derivatives, Kluwer
Academic, Boston, USA, 1991.

\bibitem{BY1} B.C. Yang, The Norm of Operator and Hilbert-Type Inequalities,
Science Press, Beijing, China, 2009.

\bibitem{BY2} B.C. Yang, A survey of the study of Hilbert-type inequalities
with parameters, Advances in Mathematics, 38(3) (2009), 257-268.

\bibitem{BY3} B.C. Yang, On the norm of an integral operator and
applications, J. Math. Anal. Appl., 321(2006),182-192.

\bibitem{JX} J.S. Xu, Hardy-Hilbert's inequalities with two parameters,
Advances in Mathematics, 36(2)(2007), 63-76.

\bibitem{BY4} B.C. Yang, On the norm of a Hilbert's type linear operator and
applications, J. Math. Anal. Appl., 325(2007), 529-541.

\bibitem{DX} D.M. Xin, A Hilbert-type integral inequality with the
homogeneous kernel of zero degree, Mathematical Theory and Applications,
30(2)(2010), 70-74.

\bibitem{BY5} B.C. Yang, A Hilbert-type integral inequality with the
homogenous kernel of degree 0, Journal of Shandong University (Natural
Science), 45(2)(2010), 103-106.

\bibitem{LB} L. Debnath, B.C. Yang, Recent developments of Hilbert-type
discrete and integral inequalities with applications, International Journal
of Mathematics and Mathematical Sciences, Volume 2012, Article ID 871845, 29
pages.

\bibitem{MY1} G.V. Milovanovi\'c, M.Th. Rassias, Some properties of a
hypergeometric function which appear in an approximation problem, Journal of
Global Optimization, 57(2013), 1173-1192.

\bibitem{MY2} M.Th. Rassias, B.C. Yang, On half-discrete Hilbert's
inequality. Applied Mathematics and Computation, 220(2.13), 75 - 93.

\bibitem{MY3} Th.M. Rassias, B.C. Yang, A multidimensional half - discrete
Hilbert - type inequality and the Riemann zeta function, Applied Mathematics
and Computation, 225(2013), 263-277.

\bibitem{MY4} M.Th. Rassias, B.C. Yang, On a multidimensional half -
discrete Hilbert - type inequality related to the hyperbolic cotangent
function. Applied Mathematics and Computation, 242(2013), 800-813.

\bibitem{MY5} M.Th. Rassias, B.C. Yang, A multidimensional Hilbert - type
integral inequality related to the Riemann zeta function, Applications of
Mathematics and Informatics in Science and Engineering (N. J. Daras, ed.),
Springer, New York, 417-433, 2014.

\bibitem{MY6} G.V. Milovanovi\'c, M.Th. Rassias (eds.), Analytic Number
Theory, Approximation Theory and Special Functions, Springer, New York, 2014.

\bibitem{BY6} B.C. Yang, A new Hilbert-type integral inequality, Soochow
Journal of Mathematics, 33(4)(2007), 849-859.

\bibitem{BY7} B.C. Yang, A Hilbert-type integral inequality with a
non-homogeneous kernel, Journal of Xiamen University (Natural Science),
48(2)(2008), 165-169.

\bibitem{BHY1} B. He, B.C. Yang, On a Hilbert-type integral inequality with
the homogeneous kernel of 0-degree and the hypergeometrc function,
Mathematics in Practice and Theory, 40(18)(2010), 203-211.

\bibitem{BY8} B.C. Yang, A new Hilbert-type integral inequality with some
parameters, Journal of Jilin University (Science Edition), 46(6)(2008),
1085-1090.

\bibitem{ZZX1} Z. Zeng, Z.T Xie, On a new Hilbert-type integral inequality
with the homogeneous kernel of degree 0 and the integral in whole plane,
Journal of Inequalities and Applications, Vol. 2010, Article ID 256796, 9
pages.

\bibitem{BY9} B.C. Yang, A reverse Hilbert-type integral inequality with
some parameters, Journal of Xinxiang University (Natural Science Edition),
27(6)(2010), 1-4.

\bibitem{WY1} A.Z. Wang, B.C. Yang, A new Hilbert-type integral inequality
in whole plane with the non-homogeneous kernel, Journal of Inequalities and
Applications, Vol. 2011, 2011: 123, doi:10.1186/1029-24X-2011-123.

\bibitem{DXY} D.M. Xin, B.C. Yang, A Hilbert-type integral inequality in
whole plane with the homogeneous kernel of degree -2, Journal of
Inequalities and Applications, Vol. 2011, Article ID 401428, 11 pages.

\bibitem{BHY2} B. He, B.C. Yang, On an inequality concerning a
non-homogeneous kernel and the hypergeometric function, Tamsul Oxford
Journal of Information and Mathematical Sciences, 27(1)(2011), 75-88.

\bibitem{BY10} B.C. Yang, A reverse Hilbert-type integral inequality with a
non-homogeneous kernel, Journal of Jilin University (Science Edition),
49(3)(2011), 437-441.

\bibitem{HWY} Q.L. Huang, S.H. Wu, B.C. Yang, Parameterized Hilbert-type
integral inequalities in the whole plane, The Scientific World Journal,
Volume 2014, Article ID 169061, 8 pages.

\bibitem{ZZX2} Z. Zhen, K. Raja Rama Gandhi, Z.T. Xie, A new Hilbert-type
inequality with the homogeneous kernel of degree -2 and with the integral,
Bulletin of Mathematical Sciences \& Applications, 3(1)(2014), 11-20.

\bibitem{WG} Z.X. Wang, D.R. Guo, Introduction to Special Functions, Science
Press, Beijing, China, 1979.

\bibitem{JK1} J.C. Kuang, Applied Inequalities. Shangdong Science and
Technology Press, Jinan, China, 2004.

\bibitem{JK2} J.C. Kuang, Introduction to Real Analysis. Hunan Educiton
Press, Changsha, China, 1996. $\ $
\end{thebibliography}
\end{document}